\documentclass[12pt]{amsart}
\usepackage{amsfonts}
\usepackage{amsmath}
\usepackage{amssymb}
\usepackage[mathcal]{eucal}
\usepackage{amscd}
\input xy
\xyoption{all}

\newcommand{\Ecal}{\mathcal{E}}

\newcommand{\Xcal}{\mathcal{X}}
\newcommand{\Ycal}{\mathcal{Y}}

\newcommand{\R}{\mathbb{R}}
\newcommand{\C}{\mathbb{C}}
\newcommand{\N}{\mathbb{N}}
\newcommand{\T}{\mathbb{T}}

\newcommand{\Z}{\mathbb{Z}}

\newcommand{\XT}{$(X,T)$\ }

\newcommand{\al}{\alpha}
\newcommand{\Ga}{\Gamma}
\newcommand{\ga}{\gamma}
\newcommand{\del}{\delta}
\newcommand{\Del}{\Delta}

\newcommand{\la}{\lambda}

\newcommand{\XGa}{$(X,\Ga)\ $}

\newcommand{\br}{\vspace{4 mm}}

\newcommand{\rest}{\upharpoonright}

\newcommand{\htop}{h_{{\rm top}}}

\newcommand{\ol}{\overline}

\newcommand{\inte}{\rm{int\,}}
\newcommand{\cls}{\rm{cls\,}}

\newcommand{\Aut}{\rm{Aut\,}}

\newcommand{\supp}{\rm{supp\,}}

\newcommand{\card}{\rm{card\,}}

\newcommand{\spann}{\rm{span}}

\newcommand{\Homeo}{\rm{Homeo\,}}

\swapnumbers
\theoremstyle{plain}
\newtheorem{thm}{Theorem}[section]
\newtheorem{cor}[thm]{Corollary}
\newtheorem{lem}[thm]{Lemma}
\newtheorem{prop}[thm]{Proposition}

\theoremstyle{definition}
\newtheorem{defn}[thm]{Definition}

\newtheorem{rmk}[thm]{Remark}

\newtheorem{prob}[thm]{Problem}

\numberwithin{equation}{section}

\begin{document}

\title
[tame dynamical systems]{On tame dynamical systems}

\author[E. Glasner]{E. Glasner}
\address{Department of Mathematics,
Tel-Aviv University, Ramat Aviv, Israel}
\email{glasner@math.tau.ac.il}

\begin{abstract}
A dynamical version of the Bourgain-Fremlin-Talagrand
dichotomy shows that the enveloping semigroup of a
dynamical system is either very large and contains
a topological copy of $\beta \N$, or it is a
``tame" topological space whose topology is determined by
the convergence of sequences. In the latter case we say
that the dynamical system is tame.
We show that (i) a metric distal minimal system is tame iff it is
equicontinuous (ii) for an abelian acting group a tame metric minimal
system is PI (hence a weakly mixing minimal
system is never tame), and (iii) a tame minimal cascade
has zero topological entropy.
We also show that for minimal distal-but-not-equicontinuous systems
the canonical map from the enveloping operator
semigroup onto the Ellis semigroup is never an isomorphism.
This answers a long standing open question.
We give a complete characterization of minimal
systems whose enveloping semigroup is metrizable.
In particular it follows that for abelian acting group
such a system is equicontinuous.
\end{abstract}

\thanks{{\em 2000 Mathematical Subject Classification
54H20}}

\keywords{Enveloping semigroup, Enveloping operator semigroup,
tame, injective dynamical system, minimal system,
equicontinuous system, Rosenthal compact,
Fr\'echet compact.}

\maketitle


\section*{Introduction}

The {\em enveloping} (or {\em Ellis }) semigroup
of a dynamical system was
introduced by R. Ellis in \cite{E1}. It proved to be
an indispensable tool in the abstract theory of topological
dynamical systems (see e.g. Ellis
\cite{E2}). However explicit computations of enveloping
semigroups are quite rare. Some examples are to be found in
Namioka \cite{Na} (1984), Milnes \cite{Mil1} (1986)
and \cite{Mil2} (1989),
Glasner \cite{Gl1} (1976) and \cite{G3}
(1993),
Berg, Gove \& Hadad \cite{BGH} (1998),
Budak, I\c sik, Milnes \& Pym. \cite{BIMP} (2001), and
Glasner \& Megrelishvili \cite{GM} (2004).
Rarely is the enveloping semigroup metrizable
(a notable exception is the case of weakly almost
periodic metric systems; see Downarowicz \cite{Dow} (1998) and Glasner
\cite{Gl2} (2003), Theorem 1.48).

In an interesting paper \cite{Ko}, A. K\"ohler pointed
out the relevance of a theorem of Bourgain, Fremlin \& Talagrand
\cite{BFT} to the study of enveloping
semigroups. She calls a dynamical system, $(X,\phi)$,
where $X$ is a compact Hausdorff space and $\phi: X \to X$
a continuous map, {\em regular} if for every function
$f\in C(X)$ the sequence $\{f\circ \phi^n: n\in \N\}$ does
not contain an $\ell^1$ sub-sequence (the sequence $\{f_n\}_{n=1}^\infty$
is an {\em $\ell^1$ sequence} if there are strictly positive constants
$a$ and $b$ such that
$$
a\sum_{k=1}^n |c_k| \le \left\| \sum_{k=1}^n c_k f_k \right\|
\le b\sum_{k=1}^n |c_k|
$$
for all $n\in \N$ and $c_1,\dots,c_n\in \C$).
Since the word ``regular" is already overused in topological dynamics
I will call such systems {\em tame}.
It turns out that for a metric system $(X,\phi)$
this is the same as the condition that $E(X,\phi)$, the enveloping
semigroup of $(X,\phi)$, be a Rosenthal compact (see \cite{GM}).

In that paper K\"ohler also considers another useful notion, that
of the {\em enveloping operator semigroup\/}.
For a Banach space $K$ and a bounded linear operator
$T: K \to K$ this is defined as
$$
\Ecal(T)={\cls}_{w^*}\{T^n :n \in \N\}.
$$
K\"ohler shows that when $(X,\phi)$ is a dynamical system,
$K= C(X)$, and $T:C(X)^* \to C(X)^*$ is the operator induced
by $\phi$ on the dual space $C(X)^*$, then there is always a
surjective homomorphism of dynamical systems
$$
\Phi: \Ecal(T) \to E(X,\phi).
$$
If we view $M(X)$, the compact space of probability
measures on $X$ equipped with the weak$^*$ topology,
as a subset of $C(X)^*$ with ${\spann}(M(X))=C(X)^*$,
we see that this map $\Phi$ is nothing
but the restriction of an element of $\Ecal(T)$ to the subspace
of Dirac measures $\{\del_x: x\in X\}$.
Theorem 5.3 of \cite{Ko} says that for a tame metric dynamical system
$(X,\phi)$, the map $\Phi$ is an isomorphism of the enveloping
operator semigroup onto the Ellis semigroup. (We will
re-prove this theorem in section 1, Theorem \ref{tame-inj}.)
In this paper I will call a dynamical system $(X,\phi)$
for which $\Phi$ is an isomorphism, an {\em injective system}.
In \cite{Ko} there are several other cases where systems are
shown to be injective and the author raises the question
whether this is always the case.
As she points out this question was posed
earlier by J. S. Pym (see \cite{Py}). In
\cite{Im} S. Immervoll gives an example of
a dynamical system which is not injective.
His example is of the form $(X,H)$ where $X=[0,1]$
is the unit interval and $H$ is an uncountable semigroup
of continuous maps from $X$ to itself.
This leaves the question open for $\Z$ (or $\N$) systems,
for group actions and for minimal systems.

In the present work the setup is that of
a compact (mostly metrizable) $\Ga$ dynamical system $(X,\Ga)$
where $\Ga$ is an arbitrary topological group.
In \cite{GM} we have shown that metrizable weakly almost periodic
(WAP) systems and more generally metrizable hereditarily almost
equicontinuous (HAE) systems are tame.
However, most of the results presented here are concerned with the case
where $(X,\Ga)$ is a {\em minimal\/} dynamical system.
In the first section it is shown that
a tame dynamical system is injective.
This, in conjunction with a theorem of Ellis and an old work
of mine on affine dynamical systems (\cite{G2}), is used to deduce
that a metric distal minimal system is injective iff
it is equicontinuous. It therefore follows that
every metric minimal distal-but-not-equicontinuous system
serves as a counterexample to the question of Pym
and K\"ohler. It is also shown that a tame minimal cascade \XT
has zero topological entropy.
In the second section I show that for abelian $\Ga$
a metric minimal tame system is PI, hence in particular
a minimal weakly mixing $\Ga$-system is never tame.
In the third section I consider the case when \XGa
is minimal and $E=E(X,\Ga)$ is metrizable.
Under these assumptions it is shown that
there is a unique minimal ideal $I$ in $E$, that
the group $K$ of automorphisms of the system $(I,\Ga)$
is compact, and that the quotient dynamical system $(I/K,\Ga)$
is proximal.
If we also assume that $\Ga$ is abelian then
$(X,\Ga)$ is equicontinuous.
In the last section I consider the question: how
big can $E(X,\Ga)$ be in $X^X$.

The reader is refered to the sources \cite{E2},
\cite{Gl1}, \cite{V}, \cite{Bro}, \cite{Au}, \cite{Vr} and
\cite{G-lu}, on the abstract theory of topological dynamics and the
structure theory of minimal dynamical systems including
the notion of PI (proximal-isometric) systems.

The questions treated in this paper arose during the work on
another one, \cite{GM}, which I just finished writing jointly
with Michael Megrelishvili. I owe him much for fruitful discussions
on these subjects. I am also indebted to Benjy Weiss
for helpful conversations; in particular the content of
section 4 was the subject of a conversation over lunch
several years ago.

\br

\section{Tame systems are injective}

Recall that a topological space $K$ is called a {\em Rosenthal compact\/}
\cite{Godefroy} if it is homeomorphic to a pointwise compact
subset of the space $B_1(X)$ of functions of the first Baire class
on a Polish space $X$. All metric compact spaces are Rosenthal. An
example of a separable non-metrizable Rosenthal compact is the
{\em Helly compact} of all
(not only strictly) increasing selfmaps
of $[0,1]$ in the pointwise topology.
Another is the {\em two arrows} space
of Alexandroff and Urysohn (see Engelking \cite{Eng}).

A topological space $K$ is a {\em Fr\'echet space} if for
every $A\subset K$ and every $x\in \ol{A}$ there exists a
sequence $x_n\in A$ with $\lim_{n\to\infty} x_n=x$
(see Engelking \cite{Eng}).
A topological space $K$ is {\em angelic}
if every relatively countably compact subset $A \subset K$
satisfies the properties (i) $A$ is relatively compact and
(ii) for every $x\in \ol{A}$
there exists a sequence $x_n\in A$ with $\lim_{n\to\infty} x_n=x$.
Thus a compact space is angelic iff it is Fr\'echet.
Clearly, $\beta \N$, the Stone-\v{C}ech compactification of the
natural numbers $\N$, cannot be embedded into a Fr\'echet space.

The following theorem is due to Bourgain, Fremlin and Talagrand
\cite[Theorem 3F]{BFT}, generalizing a result of Rosenthal. The
second assertion (BFT dichotomy) is presented as in the book of
Todor\u{c}evi\'{c} \cite{T-b} (see Proposition 1 of section 13).

\begin{thm} \label{BFT}
\begin{enumerate}
\item
Every Rosenthal compact space $K$ is angelic.
\item (BFT dichotomy)
Let $X$ be a Polish space and let $\{f_n\}_{n=1}^\infty \subset
C(X)$ be a sequence of real valued functions which is {\em
pointwise bounded}  (i.e. for each $x\in X$ the sequence
$\{f_n(x)\}_{n=1}^\infty$ is bounded in $\R$).
Let $K$ be the pointwise closure of $\{f_n\}_{n=1}^\infty$
in $\R^X$. Then either $K \subset B_1(X)$ (i.e. $K$ is
Rosenthal compact) or $K$ contains a homeomorphic copy of $\beta\N$.
\end{enumerate}
\end{thm}

The following dynamical BFT dichotomy is derived in \cite{GM}.

\begin{thm}[A dynamical BFT dichotomy]\label{D-BFT}
Let $(X,\Ga)$ be a metric dynamical system and let $E=E(X,\Ga)$
be its enveloping semigroup. We have the following alternative.
Either
\begin{enumerate}
\item
$E$ is a separable Rosenthal compact, hence with cardinality
${\card}{E} \leq 2^{\aleph_0}$; or
\item
the compact space $E$ contains a homeomorphic
copy of $\beta\N$, hence ${\card}{E} = 2^{2^{\aleph_0}}$.
\end{enumerate}
\end{thm}

\begin{defn}
We will say that an enveloping semigroup $E(X,\Ga)$ is {\em tame}
if it is separable and Fr\'echet.
A dynamical system $(X,\Ga)$ is {\em tame} when $E(X,\Ga)$ is tame.
\end{defn}

In these terms Theorem \ref{D-BFT} can be rephrased as
saying that {\em a metric dynamical system $(X,\Ga)$ is either
tame or $E(X,\Ga)$ contains a topological copy of $\beta \N$}.
When \XGa is a metrizable system
the group $\Ga$ is embedded in the Polish
group ${\Homeo}(X)$ of homeomorphisms of $X$
equipped with the topology of uniform convergence.
From this fact it is easy to deduce that the enveloping
semigroup $E(X,\Ga)$ is separable.
If moreover \XGa is tame then
$E=E(X,\Ga)$ is Fr\'echet and every element $p\in E$
is a limit of a sequence of elements of $\Ga$,
$p=\lim_{n\to \infty} \ga_n$.

Examples of tame dynamical systems include
metric minimal equicontinuous systems, almost periodic (WAP) systems
(E. Akin, J. Auslander, and K. Berg \cite{AAB}),
and hereditarily non-sensitive (HNS) systems
(Glasner and Megrelishvili \cite{GM}).

The cardinality distinction between the two cases entails
the first part of the following proposition.

\begin{prop}
\begin{enumerate}
\item
For metric dynamical systems tameness is preserved by taking
\begin{enumerate}
  \item subsystems,
  \item countable self products, and
  \item factors.
\end{enumerate}
\item
Every metric dynamical system \XGa admits a unique maximal
tame factor.
\end{enumerate}
\end{prop}

\begin{proof}
As pointed out, the first statement follows from
cardinality arguments (note that $E(X,\Ga)=E(X^\kappa,\Ga)$
for any cardinal number $\kappa$). To prove the second use Zorn's
lemma, the first part of the theorem, and the fact that
a chain of factors of a metric system is necessarily
countable, to find a maximal tame factor. Then use the first
part again to deduce that such maximal factor is unique.
\end{proof}

As was mentioned in the introduction the following
theorem is due to K\"ohler; our proof though is
different (see also \cite{Gl2}, Lemma 1.49).

\begin{thm}\label{tame-inj}
Let \XGa be a metric tame dynamical system.
Let $M(X)$ denote the compact convex
set of probability measures on $X$ (with the weak$^*$
topology). Then each element $p\in E(X,\Ga)$ defines
an element $p_*\in E(M(X),\Ga)$ and the map
$p \mapsto p_*$ is both a dynamical system and a semigroup
isomorphism of $E(X,\Ga)$ onto $E(M(X),\Ga)$.
\end{thm}

\begin{proof}
Since $E(X,\Ga)$ is Fr\'echet we have for every $p\in E$ a sequence
$\ga_i\to p$ of elements of $\Ga$ converging to
$p$. Now for every $f\in C(X)$ and every
probability measure $\nu\in M(X)$ we get by the
Riesz representation theorem and
Lebesgue's dominated convergence theorem
$$
\ga_i\nu(f)=\nu(f\circ \ga_i)\to \nu(f\circ p):=p_*\nu(f).
$$
Since the Baire class 1 function $f\circ p$ is well defined
and does not depend upon the choice of the convergent
sequence $\ga_i\to p$, this defines the map $p \mapsto p_*$ uniquely.
It is easy to see that this map is an isomorphism of dynamical systems,
whence a semigroup isomorphism. Finally as $\Ga$ is dense in both
enveloping semigroups, it follows that this isomorphism is onto.
\end{proof}

\begin{defn}
We will say that the dynamical system $(X,\Ga)$ is
{\em injective} if the natural map
$E(M(X),\Ga) \to E(X,\Ga)$ is an isomorphism.
\end{defn}

In these terms the previous theorem can be restated as
follows. {\em A tame dynamical system is injective}.
Our next theorem, which relies on \cite{G1},
answers a question of J. S. Pym
and A. K\"{o}hler (see also S. Immervoll \cite{Im}).

\begin{thm}\label{distal}
A minimal distal metric dynamical system is
injective iff it is equicontinuous.
\end{thm}

\begin{proof}
It is well known that when $(X,\Ga)$ is equicontinuous,
$E=E(X,\Ga)$ is a compact topological group and in that case
it is easy to see that $(X,\Ga)$ is injective.
By a theorem of Ellis (see e.g. \cite{E2}),
a system $(X,\Ga)$ is distal iff $E(X,\Ga)$
is a group. Thus, if $(X,\Ga)$ is distal metric and injective
then $E(X,\Ga)=E(M(X),\Ga)$ is a group and it follows that
the dynamical system $(M(X),\Ga)$ is also distal. By Theorem 1.1
of \cite{G1}, the system $(X,\Ga)$ is equicontinuous.
\end{proof}

\begin{cor}\label{cor-distal}
A  minimal distal metric system is tame iff it is equicontinuous.
\end{cor}

\begin{proof}
A metric minimal equicontinuous system is isomorphic to its
own enveloping semigroup. For the other direction observe
that if \XGa is tame then by Theorem \ref{tame-inj} it
is injective hence, by Theorem \ref{distal}, it is
equicontinuous.
\end{proof}

By way of illustration consider, given an irrational
number $\al\in \R$, the minimal distal dynamical $\Z$-system
on the two torus $(T,\T^2)$ given by:
$$
T(x,y)=(x+\al, y+x)\qquad \pmod 1.
$$
Since this system is not equicontinuous Theorem \ref{distal}
and Corollary \ref{cor-distal}
show that it is neither tame nor injective.

\br

The fact that tame systems are injective also yield the result that
metric tame minimal systems have zero topological entropy.
For this we need the following (simplified version of a) theorem of
Blanchard, Glasner, Kolyada and Maass \cite[Theorem 2.3]{BGKM}.
Recall that a pair of
points $\{x,y\} \subseteq X$ is said to be a
{\em Li--Yorke pair\/} if one has simultaneously
$$ \limsup_{n\to
\infty} d(T^nx,T^ny)=\delta >0,\quad \text {{ and }}
\quad \liminf_{n\to
\infty} d(T^nx,T^ny)=0.
$$
In particular a Li--Yorke pair is proximal.
A set $S\subseteq X$ is called {\em scrambled\/} if any pair
of distinct points $\{x,y\}\subseteq S$ is a Li--Yorke pair.
A dynamical system $(X,T)$ is called
{\em chaotic in the sense of Li and Yorke\/} if $X$ contains an
uncountable scrambled set.

\begin{thm}\label{BGKM}
Let $(X,T)$ be a topological dynamical system with $\htop (X,T)>0$. Let
$\mu$ be a $T$-ergodic probability measure with $\supp(\mu)=X$ and
$h_\mu(X,T)>0$.
Then there exists a topologically transitive subsystem $(W,T\times T)$
with $W\subseteq X\times X$, such that for every open
$U\subseteq X$ there exists a
Cantor scrambled set $K\subseteq U$ with
$K\times K\setminus \Del_X\subseteq W_{tr}$,
where $W_{tr}$ is the set of transitive points in $W$.
Thus a dynamical system with positive topological entropy
is chaotic in the sense of Li and Yorke.
\end{thm}

We note that the set $W$ in Theorem \ref{BGKM}
has the following special form.
There exists a measure theoretical weakly mixing
factor map $\pi:(X,\Xcal,\mu,T) \to (Y,\Ycal,\nu,T)$
with a corresponding measure disintegration
$$
\mu = \int_Y \mu_y\, d\nu(y)
$$
having the property that $\mu_y$ is non-atomic for $\nu$-a.e. $y$.
The subsystem $W$ is then given as $W = \supp(\la)$, where
$$
\la = \mu\underset{\nu}{\times}\mu =
\int_Y \mu_y \times \mu_y \, d\nu(y).
$$
Consequently if $X_0\subset X$ is any $\mu$-measurable set
with $\mu(X_0)=1$ then with no loss of generality we can assume that
for $\nu$ almost every $y$ the measure $\mu_y$ satisfies
the condition $\mu_y(X_0)=1$. It then follows that
the Cantor set in Theorem \ref{BGKM} can be
chosen to be a subset of $X_0$.

\begin{thm}\label{entropy}
A minimal metric tame $\Z$ dynamical system \XT has zero
topological entropy.
\end{thm}

\begin{proof}
By the variational principle it suffices to show that
$h_\mu(T)=0$ for every $T$-invariant probability measure
$\mu$ on $X$. Let $\mu$ be such a measure. By Theorem
\ref{tame-inj} \XT is injective and therefore $v_*(\mu)=\mu$
for any minimal idempotent $v\in E=E(X,T)$.
Since $v_*\mu(f)=\mu(f\circ v)$ for every $f\in C(X)$
it follows that $\mu(vX)=1$.
(Note that $vX$ is an analytic set hence universally measurable.)
Now if $h_\mu(T)>0$ then
by Theorem \ref{BGKM}, with $X_0=vX$, there is a
Cantor set $K\subset X_0$ such that for
every $x,x'$ distinct points in $K$
the pair $(x,x')$ is a proximal pair.
However, since pairs $(x,x')\in vX \times vX$
with $x\not= x'$ are {\em almost periodic\/}
(i.e. have minimal orbit closure in $X\times X$)
they are never proximal pairs and we conclude that $h_\mu(T)=0$.
\end{proof}

\begin{rmk}
In the proof of Theorem \ref{entropy}, with slight modifications,
one can use instead of the results in \cite{BGKM} a theorem
of Blanchard, Host and Ruette \cite{BHR} on the abundance of asymptotic
pairs is a system \XT with positive topological entropy.
\end{rmk}

\section{Minimal tame systems are PI}

As we have seen, when \XGa is a metrizable tame system
the enveloping semigroup $E(X,\Ga)$ is a separable Fr\'echet
space. Therefore each element $p\in E$
is a limit of a sequence of elements of $\Ga$,
$p=\lim_{n\to \infty} \ga_n$. It follows that the subset
$C(p)$ of continuity points of each $p\in E$ is a
dense $G_\del$ subset of $X$. More generally, if $Y\subset X$
is any closed subset then the set $C_Y(p)$ of continuity points
of the map $p\rest Y : Y \to X$ is a dense $G_\del$ subset
of $Y$. For an idempotent $v=v^2\in E$ we write
$C_v$ for $C_{\ol{vX}}(v)$.

\begin{lem}\label{Cp}
Let \XGa be a metrizable tame dynamical system,
$E=E(X,\Ga)$ its enveloping semigroup.
\begin{enumerate}
\item
For every $p\in E$ the set $C(p)\subset X$
is a dense $G_\del$ subset of $X$.
\item
For every minimal idempotent $v\in E$, we have
$C_v\subset vX$.
\item
When $\Ga$ is commutative we have $C(v)\subset vX$.
\end{enumerate}
\end{lem}

\begin{proof}
1.\
See the remark above.

2.\
Given $x\in C_v$ choose a sequence $x_n \in vX$ with
$\lim_{n\to\infty} x_n =x$.
We then have $vx=\lim_{n\to\infty} vx_n = \lim_{n\to\infty} x_n =x$,
hence $C_v\subset vX$.

3.\
When $\Ga$ is commutative we have $\ga p =p \ga$ for every
$\ga\in \Ga$ and $p\in E$. In particular the subset
$vX$ is $\Ga$ invariant hence dense in $X$. Thus $\ol{vX}=X$,
hence $C(v) = C_v \subset vX$ by part 2.
\end{proof}

We next proceed to the main theorem of this section.

\begin{defn}
Let \XGa be a dynamical system. We say that a closed
$\Ga$-invariant set $W\subset X\times X$ is a {\em $M$-set}
if it satisfies the conditions
\begin{enumerate}
  \item The system $(W,\Ga)$ is topologically transitive.
  \item The almost periodic points are dense in $W$.
\end{enumerate}
\end{defn}

A theorem of Bronstein asserts that a metric system \XGa
is PI iff every $M$-set in $X\times X$ is minimal
(\cite{Bro}, see also \cite{G-lu}).

\begin{thm}\label{PI}
Let $\Ga$ be a commutative group. Then
any metric tame minimal system \XGa is PI.
\end{thm}

\begin{proof}
We will prove that the Bronstein condition holds; i.e.
that every $M$-set in $X\times X$ is minimal. So let
$W\subset X\times X$ be an $M$-set.
Let $v=v^2$ be some minimal idempotent in $E(X,\Ga)$.
By Theorem \ref{Cp}.3 the set $C(v)$ of continuity points
of the map $v: X \to X$ is a dense $G_\del$ subset of $X$
and moreover $C(v)\subset vX$. Let $U$ be a relatively
open subset of $W$, then there exists a minimal subset
$M\subset W$ with $M\cap U\not =\emptyset$. Let
$\pi_i: M \to X,\ i=1,2$ denote the projection maps.
Because $M$ is minimal we have $\pi_i(M)=X$ and
the map $\pi_i$ is semi-open; i.e.
${\inte}(\pi_i(\ol{V}))\ne \emptyset$ for every nonempty open
subset $V$ of $M$ (see e.g. \cite[Lemma 1.5]{G1.5};
these observations are due to Auslander and Markley).
It follows that the sets $\pi_i^{-1}(C(v)),\ i=1,2$
are dense $G_\del$ subsets of $M$ and therefore so is
the set
$$
(C(v) \times C(v)) \cap M = \pi_1^{-1}(C(v))\cap \pi_2^{-1}(C(v)).
$$
In particular $U \cap (C(v) \times C(v)) \ne \emptyset$
and we conclude that $W_0=(C(v) \times C(v))\cap W$ is a
dense $G_\del$ subset of $W$.

Let $W_{tr}$ be the dense $G_\del$ subset
of transitive points in $W$ and observe that $W_0\cap W_{tr}
\ne\emptyset$. If $(x,x')$ is a point in $W_0\cap W_{tr}$, then
$\ol{\Ga(x,x')}=W$ and since $(x,x') \in (vX)\times (vX)$
it follows that $W$ is minimal.
\end{proof}

\begin{cor}\label{wm}
Let $\Ga$ be a commutative group and \XGa a minimal weakly
mixing metric tame dynamical system, then \XGa is trivial.
\end{cor}

\begin{proof}
A minimal system which is weakly mixing and PI is necessarily
trivial.
\end{proof}

A direct proof of Corollary \ref{wm} that does not
require the PI theory is as follows.

\begin{proof}
Fix a minimal idempotent $u\in E$ and let $C(u)\subset X$ be
the dense $G_\del$ subset of continuity points of $u$.
Fix some $x\in X$; then, by a theorem of Weiss, $P[x]$
the proximal cell of $x$,
is also a dense $G_\del$ subset of $X$
(see \cite[Theorem 1.13]{Gl2}).
Set $A=C(u)\cap P[x]$, then
for $y\in A$ there is a sequence $\ga_j \in \Ga$
such that $\lim_{j\to\infty} \ga_j x = \lim_{j\to\infty} \ga_j y
= x$. By the continuity of $u$ at $x$ we have
\begin{align*}
ux&=u\lim_{j\to\infty} \ga_j x =\lim_{j\to\infty} \ga_j ux\\
&=u\lim_{j\to\infty} \ga_j y =\lim_{j\to\infty} \ga_j uy,
\end{align*}
so that $(ux,uy)\in P$. This implies $ux=uy$ and we conclude
that $ux=uy$ for every $y\in A$. For an arbitrary element
$\ga\in \Ga$, the set $\ga^{-1} A\cap A$ is a residual subset of $X$
and for each $y$ in this set we get
$ux=u \ga y= \ga uy= uy$. Since $\Ga$ is commutative and $(\Ga,X)$ is
minimal we conclude that $\ga z =z$ for every $z\in X$. Thus
$\Ga$ acts trivially on $X$ and the minimality of $(\Ga,X)$
implies that $X$ is a one point space.
\end{proof}

In \cite{GM} there is an example of a minimal tame dynamical
cascade (i.e. a $\Z$-system) on the Cantor set with
an enveloping semigroup which is not metrizable
(in fact $E$ in this example is homeomorphic to
the ``two arrows" space). This system
has the structure of an almost 1-1 (hence proximal)
extension of an irrational rotation on the circle $\T$.
Another such example is in R. Ellis \cite{Ellis}
where the enveloping semigroup of the $SL(2,\R)$ action
on the projective line $\mathbb{P}$ is shown to be tame but not
metrizable. Here the system $(\mathbb{P},SL(2,\R))$ is proximal.
In view of these examples, Corollary \ref{cor-distal}, Theorem
\ref{entropy}, Corollary \ref{wm},
and Theorem \ref{metric} below,
it is reasonable to raise the following question.

\begin{prob}
Is it true that every minimal metrizable tame system \XGa
with an abelian acting group
is a proximal extension of an equicontinuous system?
(Or, for the general acting group, $X$ is proximally equivalent to
a factor of an isometric extension of a proximal system.)
\end{prob}

\br

\section{Metrizable enveloping semigroups}

In this section we consider the case of a minimal
dynamical system for which $E=E(X,\Ga)$ is metrizable.
Of course then $E$ is tame and if $I\subset E$
is a minimal (left) ideal in $E$ then
the dynamical system $(I,\Ga)$ is metric
with $E(I,\Ga) \cong E(X,\Ga)$ so that it is
also tame.

\begin{thm}\label{metric}
Let \XGa be a minimal dynamical system such that $E=E(X,\Ga)$ is
metrizable. Then
\begin{enumerate}
\item
There is a unique minimal ideal $I\subset E=E(X,\Ga)\cong E(I,\Ga)$,
\item
The Polish group $G_U={\Aut}(I,\Ga)$, of automorphisms
of the system $(I,\Ga)$ equipped with the topology of
uniform convergence, is compact.
\item
The quotient dynamical system $(I/K,\Ga)$ is proximal.
\item
The quotient map $\pi:I \to I/K$ is a $K$-extension.
\item
If in addition \XGa is incontractible then $I=K$ and
$\Ga$ acts on $K$ by translations via a continuous
homomorphism $J:\Ga \to K$ with $J(\Ga)$ dense in $K$.
In particular $(I,\Ga)$, and hence also \XGa, is equicontinuous.
\item
If $\Ga$ is commutative then $X=I=K$ and
$K={\cls}J(\Ga)$ is also commutative.
\end{enumerate}
\end{thm}

\begin{proof}
We split the proof into several steps.

1.\
If $I\subset E$ is a minimal left ideal then \XGa is a factor
of the dynamical system $(I,\Ga)$ and the enveloping
semigroup $E(I,\Ga)$ is isomorphic to $E(X,\Ga)$,
where each $p\in E$ is identified with the map
$L_p:E \to E,\ q\mapsto pq$.

2.\
Let $u=u^2$ be a fixed idempotent in $I$ and, as usual
denote $G=uI\subset I$. Then to each $\al \in G$
corresponds an automorphism $\hat\al: I \to I$,
which is defined by $\hat\al (p) =p\al, \forall \, p\in I$.
The map $G \to G_U,\ \al \mapsto \hat\al$ is a surjective algebraic
isomorphism. The inverse map, $G_U \to G\subset I$ is given by
$\hat\al \mapsto \hat\al(u)=u\al=\al$.
Thus $G_U$ acts on $I$ by right multiplication. In the sequel
we will identify $\hat\al$ with $\al$.

3.\
By Lemma \ref{Cp}, for each idempotent $v=v^2\in I$,
the restricted map $v: \ol{vI} \to vI$,\ $q \mapsto vq$
has a dense $G_\del$ subset $C_v\subset \ol{vI}$ of continuity points.
Again by Lemma \ref{Cp} $C_v \subset  vI$.
Since clearly $C_vG \subset C_v$, we get $C_v=vG=vI$.

4.\
If $p\in vI$ then also $p: \ol{vI} \to vI$ and thus its
set of continuity points $C_p$ is also a dense $G_\del$
subset of $\ol{vI}$. Therefore $C_p\cap C_v = C_p \cap vI\ne\emptyset$,
and since $C_pG\subset C_p$ we conclude that $C_p\supset vI$.

5.\
The $G$-dynamical system $(I,G)$ admits a minimal subset $M$,
and it is clearly of the form $M=\ol{vI}=\ol{vG}$ for some $v=v^2\in I$.
By minimality we have $M=\ol{wI}$ for any other idempotent
$w=w^2\in \ol{vI}$. Since, by step 3, $C_v=vI$ and
$C_w=wI$ are residual subsets of $M$, their intersection
is non empty and the structure of $I$ as a disjoint union
of groups implies that $v=w$, hence $wI = vI =M$.
Thus $v: I \to I, \ p\mapsto v p$
has a closed range $vI$ and the right action
of $G$ on $M=vI=vG$ is algebraically transitive.
(The right action of $G$ on $I$,
hence also on $vI$ is free.)
Moreover, from step 4 we see that every $\al\in G$
acts continuously on $vI$ on the {\em left\/}; i.e.
$p_n \to p$, for $p_n,p\in vI$ implies $\al p_n \to \al p$.
Thus in the compact group $vG$, with the topology inherited
from $I$, both left and right multiplications are continuous.
By a theorem of Ellis (\cite{E0}) it follows that $vG$ is a compact
topological group. Being a closed subset of $I$ it is
also Polish.

6.\
Now the map $v: G_U \to vG, \ \al \mapsto v\al$
is clearly a continuous surjective 1-1 homomorphism of Polish topological
groups and a theorem of Banach (\cite{B}) implies
that it is a topological isomorphism (see also \cite[Lemma 3]{G2}).
We therefore conclude that $G_U$ is a compact subgroup of
${\Aut}(I,\Ga)$. Letting $K=G_U$ all the assertions
of the theorem follow readily.
\end{proof}

\begin{rmk}
Let $\Ga$ be a topological group and $J: \Ga \to K$
a continuous homomorphism, where $K$ is a compact metrizable
topological group and $J(\Ga)$ is dense in $K$. In addition let
$H$ be a closed subgroup of $K$ for which
$\bigcap_{\, k\in K} kHk^{-1}=\{e\}$. Then the dynamical system
$(X,\Ga)=(K/H,\Ga)$ where $\ga (kH) = J(\ga)kH,\ (\ga\in \Ga, k\in K)$
is a minimal dynamical system with $E(X,\Ga)= K$.
In fact, these are the only examples I know of minimal systems
with metrizable enveloping semigroup.
\end{rmk}

\begin{prob}
Is there a nontrivial minimal proximal system
with a metrizable enveloping semigroup?
\end{prob}

\br

\section{When is $E(X,\Ga)$ all of $X^X$ ?}

We say that the system $(X,\Ga)$ is {\em $n$-complete}
if for every point $(x_1,\dots,x_n)
\in X^n$ with distinct components
the orbit  $\Ga(x_1,\dots,x_n) $ is dense in $X^n$.
It is called {\em complete\/}
when it is $n$-complete for every $n\in \N$.

\begin{thm}\label{complete}
Let $(X,\Ga)$ be a dynamical system. Then $E(X,\Ga)=X^X$
if and only if $(X,\Ga)$ is complete.
\end{thm}

\begin{proof}
Suppose $E(X,\Ga)=X^X$ and let $(x_1,\dots,x_n),
(x'_1,\dots,x'_n)\in X^n$.
Then there exists an element $p\in E$ with $px_i=x'_i,
\ i=1,\dots ,n$, hence $\ol{\Ga(x_1,\dots,x_n)} = X^n$
and $(X,\Ga)$ is complete.

Conversely if $(X,\Ga)$ is complete then
clearly every element of $X^X$ can be approximated
by an element from $\Ga$. As $E$ is closed this
concludes the proof.
\end{proof}

\begin{cor}
Let $X$ be a topological space which is $n$-homogeneous
for every $n\in \N$ (i.e. the group ${\Homeo}(X)$ acts
$n$-transitively on $X$ for every $n$) then for any
dense subgroup $\Ga\subset {\Homeo}(X)$ the dynamical system
$(X,\Ga)$ is complete, hence $E(X,\Ga)=X^X$.
For example this is the case for the Cantor set $C$, for any sphere
$S^n,\ n\ge 2$, and for the Hilbert cube $Q$.
\end{cor}

If $\phi$ is a nontrivial continuous automorphism of a
system \XGa then $\phi p=p\phi$ for every $p\in E=E(X,\Ga)$.
Thus when the group ${\Aut}(X,\Ga)$ is nontrivial
then $E\subset \{p\in X^X: \phi p=p\phi,\ \forall \phi\in
{\Aut}(X,\Ga)\}$. In particular, when $\Ga$ is commutative
$$
E\subset \{p\in X^X: p \ga=\ga p,\ \forall \ga\in
\Ga\}.
$$
Are there dynamical systems \XGa for which this
inclusion is an equality?
A $\Z$-dynamical system $(X,T)$ is said to have {\em
$2$-fold topological minimal
self-joinings\/} (\cite{J}, \cite{K}) if it satisfies
the following condition.
For every pair $(x,x')\in X \times X$ with $x'\not\in \{T^nx: n\in \Z\}$,
the orbit $\{T^n(x,x'): n\in \Z\}$ is dense in $X\times X$.
If it satisfies the analogous condition for every point
$(x_1,x_2,\dots,x_n)\in X^n$ whose coordinates belong to
$n$ distinct orbits, then \XGa has
{\em $n$-fold topological minimal self-joinings}.
As in the proof of Theorem \ref{complete} it is easy to see that
$$
E(X,T)=\{p\in X^X: p\, T=T p\}
$$
iff \XGa has {\em $n$-fold topological minimal self-joinings\/}
for all $n\ge 1$.
Now in \cite{K}, J. King shows that no non-trivial map has
$4$-fold topological minimal self-joinings. We thus get
the following.

\begin{thm}
There does not exist an infinite minimal cascade \XT for which
$$
E(X,T)=\{p\in X^X: p\, T=T p\}.
$$
\end{thm}

\begin{rmk}
In \cite{W} B. Weiss shows that every aperiodic ergodic zero entropy
measure preserving system has a topological model
which has two-fold topological minimal self-joinings
({\em doubly minimal\/} in his terminology).
\end{rmk}


\bibliographystyle{amsplain}

\end{document}